\documentclass[11pt]{article}
\usepackage[english]{babel}
\usepackage{graphics}
\usepackage{amsfonts}
\usepackage{amsmath}
\usepackage{amssymb}
\usepackage{mathrsfs}
\usepackage{enumerate}
\usepackage{hyperref}
\usepackage{epsf}
\usepackage{epsfig}
\usepackage{amscd}
\usepackage{amsmath}

\newtheorem{lemma}{Lemma}[section]
\newtheorem{corollary}{Corollary}[section]
\newtheorem{proposition}{Proposition}[section]
\newtheorem{example}{Example}[section]
\newtheorem{theorem}{Theorem}[section]

\newtheorem{remark}{Remark}[section]
\newtheorem{definition}{Definition}[section]
\def\psn{\par\smallskip\noindent}
\def\QED{\hfill $\Box$\par\smallskip\noindent}

\def\scatola{\lower5pt\hbox{\vbox{\hrule\hbox{\vrule\kern2pt\vbox%
{\kern5pt\hbox{\mathsurround=0pt }\kern2pt}\kern4pt\vrule}\hrule}}\
} 

\def\diam{{\rm diam}}

\def\H{{\mathbb H}}
\def\Hn{{\mathbb H^n}}

\def\R{{\mathbb{R}}}
\def\G{\mathbf{G}}

\def\dom{\texttt{\rm dom}}
\def\gr{\texttt{\rm gr}}

\def\0{\overline 0}

\def\psn{\par\smallskip\noindent}

\font\tenmsb=msbm10 \font\sevenmsb=msbm7 \font\fivemsb=msbm5
\newfam\msbfam
\textfont\msbfam=\tenmsb \scriptfont\msbfam=\sevenmsb
\scriptscriptfont\msbfam=\fivemsb

\begin{document}

\title{On the local boundedness of maximal H--monotone operators}
\author{Z.M. Balogh\thanks{Institute of Mathematics, University of Bern, Sidlerstrasse 5,
3012 Bern, Switzerland {(\tt zoltan.balogh@math.unibe.ch)}, ZMB was partially supported by
the Research Grant: Geometric Analysis of sub-Riemannian Spaces, Proj. Nr. 200020-165507 of the  Swiss National Science Foundation.},
 A. Calogero\thanks{
Dipartimento di Matematica e Applicazioni, Universit\`a degli Studi
di Milano--Bicocca, Via Cozzi 55, 20125 Milano, Italy ({\tt
andrea.calogero@unimib.it}, {\tt rita.pini@unimib.it})},
 R. Pini$^\dag$
 }
\date{}
\maketitle

\begin{abstract}
\noindent  In this paper we prove that maximal H-monotone operators
$T:\H^n\rightrightarrows V_1$ whose domain is all the Heisenberg
group $\H^n$ are locally bounded. This implies that they are upper
semicontinuous. As a consequence, maximal H-monotonicity of an
operator on $\H^n$ can be characterized  by a suitable version of
Minty's type theorem.
\end{abstract}


\medskip




\section{Introduction}
Maximal monotone maps in Euclidean spaces $\R^n$ and, more in
general, in Hilbert spaces, play key roles in  evolution equations
and in other fields of functional analysis. The most notable example
of a maximal monotone map in $\R^n$ is provided by the
subdifferential map $\partial f$ associated to a convex function
$f:\R^n\to\R.$

The celebrated Minty theorem provides a characterization of maximal
monotonicity (see \cite{RoWe2004}): given a monotone set-valued map
$T:\R^n\rightrightarrows \R^n,$ then $T$ is maximal monotone if and
only if $I+\lambda T$ is surjective onto $\R^n,$ for every $\lambda
>0;$ in this case, the resolvent map $(I+\lambda T)^{-1}$ is single-valued and 1-Lipschitz continuous on $\R^n.$

For operators defined on Carnot groups $\G$, a notion of
H-monotonicity, and maximal H-monotonicity, has been introduced in
\cite{CaPi2012}. This notion fits the monotonicity of maps in
Euclidean spaces to the horizontal structure $V_1$ of $\G.$ It
arises naturally as the property fulfilled by the H-normal map
$\partial_Hf$ associated to an H-convex function $f:\G\to \R.$

In the classical case, well-known regularity properties enjoyed by
maximal monotone maps $T:\R^n\rightrightarrows\R^n$ are upper
semicontinuity and local boundedness in the interior of the domain
of $T$; in particular, the proof of the latter relies essentially on
the fact that any given ball of $\R^n$ is contained in the convex
hull of at most $n+1$ points.

In this paper, we investigate maximal H-monotone operators
$T:\H^n\rightrightarrows V_1$ defined on the Heisenberg group
$\H^n$, where $V_1\cong \R^{2n}$ denotes the first layer of the Lie
algebra of $\H^n.$ An important example of a such operator is the
horizontal normal map $\partial_H f$ of an H-convex function
$f:\H^n\to\R.$ When dealing with these operators, one has to face a
much more intricate situation, due to the lack of the Euclidean
geometry of the underlying setting. More precisely, we say that
$T:\H^n\rightrightarrows V_1$ is H-monotone if for every
$\eta\in\Hn,\ \eta' \in H_\eta,\ v\in T_\lambda(\eta)$ and $v'\in
T_\lambda(\eta')$ we have (see Definition \ref{def h mon})
$$\langle v-v',\xi_1(\eta)-\xi_1(\eta')\rangle\ge 0,$$
where $H_{\eta}$ is the horizontal plane through $\eta$ and $\xi$ is the canonical projection $\xi_1(x,y,t)=(x,y)\in V_1$,  for every
$(x,y,t)\in\H^n.$
 The restriction $\eta' \in H_\eta$ is an essential one, it implies that
the notion of H-monotonicity provides information on the behaviour of
the operator $T$ at the point $\eta\in\Hn$ only along the horizontal
directions through $\eta$. This restriction creates major difficulties in studying the
properties of H-monotone maps. Despite the fact that several
notions of convex hulls have been introduced in $\Hn$ (\cite{CaCaPi}), they seem not to be
useful for our purpose.

 The goal of this paper is to overcome the above indicated difficulties and
study  the local boundedness of maximal H-monotone maps. In Theorem
\ref{teo usc} we show that, for an operator $T$ with $\dom(T)=\H^n,$
upper semicontinuity is equivalent to local boundedness.  Our main result is the following:

\begin{theorem}\label{teo bounded}
Let $T: \H^n \rightrightarrows V_1$ be a maximal $H$-monotone map,
such that $\mathrm{dom}(T)=\H^n.$ Then $T$ is locally bounded.
\end{theorem}
This statement implies that $T$ is upper semicontinuous. The proof of
this theorem is considerably more involved when compared to the
Euclidean framework.  The statement  recovers the same result as in the Euclidean case with
considerably reduced assumptions as we can use information provided by the monotonicity
only along horizontal directions. Our proofs require  a deeper understanding of the
horizontal geometry of $\H^n$; in particular, the non-integrability of the horizontal bundle, or
the so-called {\it twirling effect} (see \cite{BaCaKr2013})
 of horizontal planes, is used repeatedly in our considerations.

Theorem \ref{teo bounded} sheds a new light on the regularity
properties of a maximal H-monotone operator on $\H^n$ and leads to
the proof that any maximal H-monotone operator on $\Hn$ can be
characterized by a suitable version of Minty's type theorem, thereby
improving a previous result by two of the authors
\cite{CaPi2014}.

\begin{theorem}\label{Minty Hn}
Let $T:\H^n\rightrightarrows V_1$ be an H-monotone map with
$\dom(T)=\H.$ Then the following two properties are equivalent:
\begin{itemize}
\item[i.]  $T$ is maximal H-monotone;

\item[ii.] for every fixed $\eta\in\Hn$ and $\lambda>0,$
the map $\left(\xi_1+\lambda T\right)\bigr|_{H_\eta}$ is surjective
onto $V_1$.
\end{itemize}
\end{theorem}
As we will see in subsection \ref{subsection resolvent}, another
application of our main theorem is the study  of the regularity of
the resolvent
 $\left(\xi_1+\lambda T\right)^{-1}:V_1 \rightrightarrows \H^n$.
 Another forthcoming application (see \cite{BaCaPiPe2016}),
 following the line of  investigation in \cite{AlAm1999}, \cite{AlAmCa}, will target
  the study of the Hausdorff dimension of singular
 sets
$\Sigma^k(T)=\{\eta\in\H^n:\ \texttt{\rm dim}(T(\eta))\ge k \}$, for
H-monotone maps $T$ and integers $k$.

\section{Basic notions and preliminary results}

\subsection{The Heisenberg group $\H^n.$}

The Heisenberg group $\Hn$ is the simplest Carnot group of step 2.
In this section we will recall some of the necessary notation and
background results used in the sequel. We will focus only on those geometric
aspects that are relevant to our paper.  For a general overview of the subject we
refer to \cite{BoLaUg2007}.

The  Lie algebra  $\frak{h}$ of $\mathbb H^n$ admits a
stratification  $ \frak{h}=V_1\oplus V_2$ with $V_1=\texttt{\rm
span}\{X_{i},\, Y_i;\ 1\le i\le n\}$ being the first layer of the so
called horizontal vector fields, and $V_2=\texttt{\rm span}\{T\}$
being the second layer which is one-dimensional. We assume
$[X_i,Y_i]=-4T$ and the remaining commutators of basis vectors
vanish. The exponential map $\exp:\frak{h}\to \mathbb H^n$ is
defined in  the usual way. By these commutator rules we obtain,
using the Baker-Campbell-Hausdorff formula, that $\Hn$ can be identified with
$\R^{n}\times \R^{n}\times \R$  endowed
with the non-commutative group law given by
$$
\eta\circ \eta'=(x,y,t)\circ (x',y',t')=(x+x',y+y',t+t'+2 (\langle
x',y\rangle-\langle x,y'\rangle)),
$$
where $x,y,x'$ and $y'$ are in $\R^n$, $t\in \R$, and for $z,z'\in
\R^n,$ we have $\langle z, z'\rangle=\sum_{j=1}^n z_j {z_j'}$
 the inner product in $\R^n$. Let us denote by $e$ the neutral element in $\Hn.$
Transporting the basis vectors of $V_1$ from the origin to an
arbitrary point of the group by a left-translation, we obtain a
system of left-invariant vector fields written as first order
differential operators as follows
\begin{equation} \label{vector fields} \left.
  \begin{array}{lll}
     X_j=\partial_{x_j}+2y_j \partial_t,\qquad j=1,...,n,\\
     Y_j=\partial_{y_j}-2x_j\partial_t,\qquad j=1,...,n.
  \end{array}
\right.
\end{equation}

Via the exponential map $\exp:\mathfrak{h}\to\H$ we identify the
vector $\sum_{i=1}^n(\alpha_i X_i+\beta_i Y_i)+\gamma T$ in
$\mathfrak{h}$ with the point $(\alpha_1,\ldots,\alpha_n,
\beta_1,\ldots,\beta_n, \gamma)$ in $\Hn;$ the inverse $\xi:
\Hn\to\mathfrak{h}$ of the exponential map has the unique
decomposition $\xi=(\xi_1,\xi_2),$ with $\xi_i:\Hn\to V_i.$ Since we
identify $V_1$ with $\R^{2n}$ when needed, $\xi_1:\Hn\to V_1\cong
\R^{2n}$ is given by $\xi_1(x,y,t)=(x,y).$

Let $N(x,y,t)=((\|x\|^2+\|y\|^2)^2+t^2)^\frac{1}{4}$ be the gauge
norm in $\mathbb H^n$. It is an interesting exercise (see \cite{Cy}) to check that
the expression
$$d_{H}(\eta,\eta')=N((\eta')^{-1}\circ \eta)$$
satisfies the triangle inequality defining a metric on $\Hn$: this
metric is the so-called Kor\'anyi-Cygan metric which is by
left-translation and dilation invariance bi-Lipschitz equivalent to
the Carnot-Carath\'eodory metric. Here, the non-isotropic Heisenberg
dilations $\delta_{\lambda}: \H^{n}\to \H^{n}$ for $\lambda > 0$ are
defined by $\delta_{\lambda}(x,y,t) = (\lambda x,\lambda y,
\lambda^{2}t)$. The Kor\'anyi-Cygan  ball of center $\eta_0\in
\mathbb H^n$ and radius $r>0$ is given by
$B_{\Hn}(\eta_0,r)=\{\eta\in \Hn:\ d_H(\eta_0,\eta)\le r\}.$

 The horizontal structure relies on the notion of horizontal
plane: given a point $\eta_0\in\Hn$, the {\emph{horizontal plane}}
$H_{\eta_0}$ associated to $\eta_0=(x_0,y_0,t_0)$ is the plane in
$\Hn$ defined by
$$
H_{\eta_0}=\left\{\eta=(x,y,t)\in\Hn:\ t=t_0+2(\langle y_0,x\rangle
-\langle x_0,y\rangle)\right\}.
$$
This is the plane spanned by the horizontal vector fields $\{X_i,\
Y_i\}_i$ at the point $\eta_0.$ We note that  $\eta'\in H_{\eta}$ if
and only if $\eta\in H_{\eta'}.$

\subsection{Multivalued maps on $\Hn$.}

Let us consider a set-valued map $T:\Hn\rightrightarrows V_1;$ we
denote by $\dom (T)$ the \emph{effective domain of $T$}, i.e. the
set $\{\eta\in\H^n:\ T(\eta)\not=\emptyset\},$ and by  $\gr(T)$ the
\emph{graph of $T$}, i.e. $\{(\eta,v)\in\H^n\times V_1:
\eta\in\dom(T), v\in T(\eta)\}.$

Let $T:\H^n\rightrightarrows V_1$ be a set-valued map, with closed
values, i.e. $T(\eta)$ is a closed set for every $\eta$. We recall
(see \cite{AlBo1999} for this general setting) that $T$ is
\emph{upper semicontinuous} (briefly \emph{usc}) at $\eta\in \H^n$
if, for every positive $\epsilon,$ there exists $\delta>0$ such that
$$
T(\eta')\subseteq T(\eta)+B_{\R^{2n}}(0,\epsilon),\quad \forall
\eta'\in \H^n,\; d_{H}(\eta',\eta)<\delta,
$$
where $T(\eta)+B_{\R^{2n}}(0,\epsilon)$ denotes the Minkowski sum of
the two sets in $\R^{2n}$. If the operator $T$ is compact-valued, i.e. $T(\eta)$
is a compact for every $\eta$, then the usc of $T$ can be
equivalently given as follows: if $\eta_k\to \eta,$ and $v_k\in
T(\eta_k),$ then there exists a subsequence $\{v_{k_n}\}$ such that
$v_{k_n}\to v\in T(\eta).$ We say that $T$ is \emph{closed} if
$\gr(T)$ is a closed subset of $\H^n\times V_1.$

Note that there is a gap between the dimension of the source and target spaces in this
definition, unlike in the Euclidean case. Nevertheless, some basic properties  follow in the same way as in
the Euclidean setting. First, the properties of being upper
semicontinuous, or closed, are related. Indeed,
\begin{remark}(see \cite{AlBo1999}, Th. 16.12)
Let $T:\H^n\rightrightarrows V_1$.  Then the following statements
hold:
\begin{enumerate}[i.]
\item if $T$ is usc and closed-valued, then it is closed;
\item if $T$ is closed, and $\mathrm{rge}(T)$ is compact, then $T$ is upper semicontinuous.
\end{enumerate}
\end{remark}

\noindent
Single-valued continuous functions map compact sets to compact sets. This property is also
true for upper semicontinuous compact-set valued maps:

\begin{proposition}(see \cite{AlBo1999}, Lemma 17.8)\label{pro bounded}
Let $T:\H^n\rightrightarrows V_1$ be a compact-valued usc map. Then
$T(K)\subset V_1$ is compact for every compact set $K\subset\Hn.$
\end{proposition}

\subsubsection{H--monotone and H--cyclical monotone maps.}

 We say that $A\subset \H^n\times V_1$ is
\emph{H--monotone} (see \cite{CaPi2014}) if
\begin{equation}\label{def set h mon}
\langle \xi_1(\eta)-\xi_1(\eta'),v-v'\rangle\ge 0,\qquad\forall
(\eta,v)\in A,\ (\eta',v')\in A,\ \eta'\in H_\eta . \end{equation}
 We stress that in the previous definition, for every point $(\xi,v)$ in the set $A$,
 the H-monotonicity property gives us  information about $A$ {\it only} in the horizontal directions
$\{X_i(\xi), Y_i(\xi)\}_i$
  through $\xi$; more precisely, (\ref{def set h mon}) is equivalent to
$$
\langle \xi_1(\eta)-\xi_1(\eta\circ\exp(tw)),v-v'\rangle\ge
0,\qquad\forall (\eta,v)\in A,\ (\eta\circ\exp(tw),v')\in A,\
t\in\R,\ w\in V_1,$$ where, for every $w$ fixed, $t\mapsto
\eta\circ\exp(tw)$ is the so called horizontal segment. This restriction gives
rise to the most difficulties of our study.

We say that $A$ is \emph{maximal H--monotone} if there are no
H--monotone sets $B\subset \H^n\times V_1$  such that $A\subset B$
and there exists $(\eta,v)\in B$ such that $ (\eta,v)\not\in A.$ As
usual, such notions of monotonicity and maximality are inherited by
the functions as follows:
\begin{definition}\label{def h mon}
We say that a set-valued map $T:\H^n\rightrightarrows V_1$ is an
H-monotone map
 if $\gr(T)$ is an H-monotone set, i.e.
for every $\eta\in\Hn,\ \eta'\in H_\eta,\ v\in T(\eta)$ and $v'\in
T(\eta')$ we have
\begin{equation}\label{def h monotone}
\langle v-v',\xi_1(\eta)-\xi_1(\eta')\rangle\ge 0.
\end{equation}
We say that $T$ is strictly H--monotone, if for every $\eta\in\H^n,\
\eta'\in H_\eta$ with $\eta'\not=\eta,\ v\in T(\eta)$ and $v'\in
T(\eta')$ in (\ref{def h monotone}) we have a strict inequality.
Moreover, we say that $T$ is maximal H--monotone
 if the set $\gr(T)$ is  maximal H--monotone.
\end{definition}

A stronger version of the concept of monotonicity is the notion of
cyclical monotonicity: in our context we say that $A\subset
\H^n\times V_1$ is an \emph{H--cyclically monotone} set (see
Definition 6.1 in \cite{CaPi2012})
 if for every sequence
$\{(\eta_i,v_i)\}_{i=0}^{m}\subset A$ such that
$\{\eta_i\}_{i=0}^{m}$ is a closed H-sequence, i.e. $\eta_i\in
H_{\eta_{i+1}}$ and $\eta_{m+1}=\eta_0,$  we have that
\begin{equation}\label{def Hcyclically}
\sum_{i=0}^{m} \langle \xi_1(\eta_{i+1}),v_i\rangle \le
\sum_{i=0}^{m} \langle \xi_1(\eta_{i}),v_i\rangle.
\end{equation}
 Moreover, we say that $A$ is \emph{maximal
H--cyclically monotone} if there are no H--cyclically monotone sets
$B\subset \H^n\times V_1$ such that $A\subset B$ and there exists
$(\eta,v)\in B$ such that $ (\eta,v)\not\in A.$ A set-valued map
$T:\H^n\rightrightarrows V_1$ is a (maximal) H-cyclically monotone
map
 if $\gr(T)$ is a (maximal) H-cyclically monotone set.

Given a function $u:\Hn\to\R$ we define the horizontal normal map of
$u$, $\partial_H u:\H^n\rightrightarrows V_1,$  by
$$
\partial_H u(\eta)=\{p\in V_1:\ u(\eta')\ge u(\eta)+\langle p,\xi_1(\eta')-\xi_1(\eta)\rangle,\ \forall\eta'\in H_\eta\}.
$$
It is well known that a function $u:\H^n\to\R$ is H--convex (see
\cite{DaGaNh2003}) if and only if $\partial_H u(\eta)$ is non empty,
for every $\eta.$ Moreover, for an H-convex function $u,$  we have
that $\partial_H u$ is H--cyclically monotone.

\noindent A cyclically monotone map has a better regularity since
essentially it coincides with the horizontal normal map of an
H-convex function. More precisely, in \cite{CaPi2012} the authors
proved that if $T:\H^n\rightrightarrows V_1$ is an H-cyclically
monotone map with $\mathrm{dom}(T)=\H^n, $ then there exists an
H-convex function $u:\H\to \R $ such that $ \gr(T)\subset
\gr(\partial_H u);$ if, in addition, $T$ is maximal, then
$\gr(T)=\gr(\partial_H u).$

 We have the following result (see \cite{CaPi2014}) of Minty type
in the case $n=1$:

\begin{theorem}\label{Minty 1}
Let $T:\H\rightrightarrows V_1$ be an H-monotone map with
$\dom(T)=\H.$
\begin{itemize}
\item[i.] If $T$ is maximal H-cyclically monotone, then the map
$(\xi_1+\lambda T)|_{H_\eta}$ is surjective onto $V_1$ for every
$\eta\in \H$ and $\lambda>0$.

\item[ii.] If the map $(\xi_1+\lambda T)|_{H_\eta}$ is
surjective onto $V_1$ for every $\eta\in \H$ for some $\lambda>0$,
then $T$ is maximal H--monotone.
\end{itemize}
\end{theorem}

\noindent Theorem \ref{Minty Hn} is a generalisation of  Theorem
\ref{Minty 1}, since we  remove the H-cyclically monotone
assumption in i., and show that the result holds in $\H^n.$
We note here, that every H--cyclically monotone set/map is an
H-monotone set/map: the following example will convince the reader
that the contrary is false, i.e. there exist maps that satisfies the
assumption in Theorem \ref{Minty Hn}, but not the assumption i. in
Theorem \ref{Minty 1}:

\begin{example}\label{example1}
Let us consider $T:\H^1\rightrightarrows V_1$ defined by
$$
T(x,y,t)=\tilde T (x,y)=(3x,-2x+4y).
$$
Then it follows (see Example 1 in \cite{CaPi2014} for the details) that $T$
is maximal H-monotone, but not maximal H-cyclically monotone.
\end{example}

\subsubsection{Usc and local boundedness for maximal H--monotone
maps.}

The purpose of this section is to establish the equivalence of usc and the local boundedness of maximal H-monotone maps.
Let us start with the following preliminary result.

\begin{proposition}\label{remark 1}
Let $T$ be a maximal H--monotone operator; then

\begin{itemize}
\item[i.] $T(\eta)$ is closed and convex (possibly empty) for every
$\eta\in\H^n;$
\item[ii.] if, in addition, $\dom(T)=\H^n,$ then
$T$ is compact-valued.
\end{itemize}
\end{proposition}

\noindent \textbf{Proof}: Let $\{v_k\}_k\subset T(\eta),$ with
$v_k\to v;$ then
$$
\langle \xi_1(\eta)-\xi_1(\eta'),v_k-v'\rangle\ge 0,\qquad\forall
\eta'\in H_\eta,\ v'\in T(\eta').$$ Taking the limit as $k\to
\infty$, we obtain
$$
\langle \xi_1(\eta)-\xi_1(\eta'),v-v'\rangle\ge 0,\qquad\forall
\eta'\in H_\eta,\ v'\in T(\eta'),$$ and the maximality implies that
$v\in T(\eta).$ This proves the closedness of $T(\eta)$. To show the
convexity, consider $v_1$ and $v_2$ in $T(\eta)$ and $\lambda\in
(0,1);$ clearly
$$
\langle \xi_1(\eta)-\xi_1(\eta'),\lambda
v_1+(1-\lambda)v_2-v'\rangle =\lambda\langle
\xi_1(\eta)-\xi_1(\eta'), v_1-v'\rangle +(1-\lambda) \langle
\xi_1(\eta)-\xi_1(\eta'),v_2-v'\rangle \ge 0,$$ for every $ \eta'\in
H_\eta,\ v'\in T(\eta').$ Again the maximality of $T$ implies that
$\lambda v_1+(1-\lambda)v_2\in T(\eta).$ Hence the proof of i. is
finished.

Let us prove ii.
 Fix $\eta\in\H^n.$ We know that $T(\eta)$ is closed: we have to show that $T(\eta)$ is bounded. Assuming the contrary,
 let us suppose that
there exists $\{v_k\}\subset T(\eta),$ such that $\|v_k\| \to
+\infty.$ Since $\{v_k\}\subset V_1,$ there exists $w\in V_1$ and a
subsequence $\{v_{k_m}\}$ such that $\langle w,v_{k_m}\rangle\to
+\infty.$ Considering the point $\eta\circ \exp w\in H_\eta,$ and any
$v\in T(\eta\circ \exp w)$  we obtain that $\langle
w, v-v_{k_m}\rangle \to -\infty,$ contradicting the H-monotonicity
of $T$.
 \QED

\noindent In particular, from the previous proposition, and from
Proposition \ref{pro bounded}, we immediately get that
\begin{corollary}\label{cor_usc_closed}
Let $T:\H^n\rightrightarrows V_1$ be a usc maximal H-monotone
operator with $\dom(T)=\H^n$. Then $T$ is closed and maps compact
sets into compact sets. In particular, it is locally
bounded.\end{corollary}

As a converse to the above Corollary, we will show that, under suitable
assumptions, local boundedness implies upper
semicontinuity. Let us first state the following technical lemma:
\begin{lemma}\label{lemma}
Let us consider $\eta, \eta'\in \H^n$ with  $\eta\not=\eta'$ and
$\eta'\in H_\eta,$ and a sequence $\{\eta_k\}_k\subset\H^n$ with
$\eta_k\to\eta$ and $\eta'\not\in H_{\eta_k}.$ Then there exists a
sequence $\{\eta'_k\}_k\subset \H^n$ with the following properties:
\begin{itemize}
\item[{a.}]
$\eta'_k\in H_{\eta'}\cap H_{\eta_k};$

\item[{b.}] $\eta'_k\to\eta';$

\item[{c.}]
$\displaystyle\frac{\xi_1(\eta'_k)-\xi_1(\eta')}{\|\xi_1(\eta'_k)-\xi_1(\eta')\|}\to
\frac{(\xi_1(\eta)-\xi_1(\eta'))}{\|\xi_1(\eta)-\xi_1(\eta')\|}.$

\end{itemize}
\end{lemma}
\noindent \textbf{Proof}:
Let us suppose, without loss of generality, that
$$
\eta=e:=(0,0,0),\qquad \eta'=(x',y',0)\not=\eta;
$$
moreover, $\eta_k=(x_k,y_k,t_k).$ Since $\eta'\in H_{e}$ and
$\eta_k\to e,$ we have that $\xi_1(\eta')\neq (0,0);$ we will
suppose that $x'\neq 0.$ Moreover, $\xi_1(\eta_k)\not=\xi_1(\eta'),$
for large $k;$ hence $H_{\eta'}\cap H_{\eta_k}\not=\emptyset.$ In
addition, $\eta_k\notin H_{\eta'},$ therefore,
$$
t_k+2(\langle x',y_k\rangle-\langle y',x_k\rangle)\neq 0.
$$
Our aim is to construct a sequence $\eta'_k$ satisfying conditions
$a.,\ b.$ and $c.$  Set $\eta'_k=(x'_k,y'_k, t'_k),$ where
\begin{equation}\label{sequence}
x'_k=(1+\epsilon_k)x',\quad y'_k=(1+\epsilon_k)y'+A_k\epsilon_k^2
x', \quad A_k=-\mathrm{sgn}(t_k+2(\langle x',y_k\rangle -\langle
y',x_k\rangle)).
\end{equation}
We will show that there exists a sequence $\{\epsilon_k\}_k$, with
$\epsilon_k>0$ and $\epsilon_k\to 0,$ such that $a-c.$ hold.
 Indeed, for such sequence
$\{\epsilon_k\}_k,$ condition $c.$ is satisfied; indeed,
$$
\frac{(\epsilon_k x', \epsilon_k y'+A_k\epsilon_k^2
x')}{\|(\epsilon_k x', \epsilon_k y'+A_k\epsilon_k^2 x')\|}=
\frac{(x',  y'+A_k\epsilon_k x')}{\|(x', y'+A_k\epsilon_k
x')\|}\to \frac{(x',  y')}{\|(x', y')\|}.
$$

Let us show that such a sequence does exist. The condition
$\eta'_k\in H_{\eta'}\cap H_{\eta_k}$ is equivalent to the
following:
\begin{equation}\label{condition_etaprimeenne}
t'_k=2(\langle y',x'_k\rangle -\langle x',y'_k\rangle
)=t_k+2(\langle y_k,x'_k\rangle -\langle x_k,y'_k\rangle ).
\end{equation}
Taking into account (\ref{sequence}), the second equality in (\ref{condition_etaprimeenne}) becomes
$$
a_kA_k\epsilon_k^2+b_k\epsilon_k+c_k=0,
$$
where
$$
a_k=(\|x'\|^2-\langle x',x_k\rangle ),\quad b_k=(\langle
x',y_k\rangle -\langle y',x_k\rangle ),\quad c_k=(t_k/2+\langle
x',y_k\rangle -\langle y',x_k\rangle ).
$$
 For every $k,$ sufficiently large, $a_k>0;$ moreover $c_k\not=0$ since
$\eta_k\not\in H_{\eta'}.$ Hence we have two solutions
$$
\epsilon_{k,\pm}=\frac{-b_k \pm\sqrt{b^2_k +4a_k|c_k|}}{2A_k a_k}.
$$
Since $c_k\to 0,$ we have  $\epsilon_{k,\pm}\to 0.$ For every $k$,
we define
$$
\epsilon_k=\begin{cases}
\epsilon_{k,+}& \mathrm{if} A_k>0\\
\epsilon_{k,-}& \mathrm{if} A_k<0
\end{cases}
$$
The sequence $\{\epsilon_k\}$ satisfies the condition $\epsilon_k>0$
and $\epsilon_k\to 0,$ therefore the sequence $\{\eta'_k\}$ defined
in \eqref{sequence} proves the assertion. \QED

\begin{theorem}\label{teo usc}
Let $T:\H^n\rightrightarrows V_1$ be maximal H-monotone, with $\mathrm{dom}(T)=\H^n.$ Then $T$ is locally bounded if and only if $T$
is usc.
\end{theorem}
\noindent \textbf{Proof}: By Corollary \ref{cor_usc_closed} we need
to prove only the \lq\lq if \rq\rq\ part. We argue by contradiction.
Suppose that $T$ is not usc. Then there exists
$\{(\eta_k,v_k)\}_k\subset\H^n\times V_1$ with
$(\eta_k,v_k)\to(\eta,v)$ with $v_k\in T(\eta_k),$ but with
$v\not\in T(\eta).$ Since $T$ is maximal, there exists a point
$\eta'\in H_\eta$ and exists $v'\in T(\eta')$ such that
\begin{equation}\label{dim pro 1}
\langle \xi_1(\eta)-\xi_1(\eta'),v-v'\rangle<0.
\end{equation}
Suppose now that there is a subsequence $\{\eta_{k_j}\}_j$ of
$\{\eta_{k}\}_k$ such that $\eta_{k_j}\in H_{\eta'}:$ then
$$
\langle \xi_1(\eta_{k_j})-\xi_1(\eta'),v_{k_j}-v'\rangle \ge0,\qquad \forall
j;
$$
taking the limit, we obtain $ \langle
\xi_1(\eta)-\xi_1(\eta'),v-v'\rangle \ge0$ which contradicts
(\ref{dim pro 1}). Hence, for large $k_0,$ $\eta_k \not\in
H_{\eta'}$ for $k\ge k_0.$ In particular
$$
\eta'\not\in H_{\eta_k},\qquad \forall k\ge k_0.
$$

Now we define the sequence $\{\eta'_k\}_k\subset \H^n$ as in Lemma
\ref{lemma}. By the local boundedness of $T,$ up to considering a
subsequence, there exists $\{v'_k\}_k,$ with $v'_k\in T(\eta'_k)$
and $v'_k\to v''\in V_1.$ Since $\eta'_k\in H_{\eta_k},$ by the
H-monotonicity of $T$ we have
$$
\langle\xi_1(\eta_k)-\xi_1(\eta'_k),v_k-v'_k\rangle\ge 0, \qquad
\forall k\ge k_0;
$$
passing to the limit we obtain
\begin{equation}\label{dim pro 2}
\langle\xi_1(\eta)-\xi_1(\eta'),v-v''\rangle\ge 0.
\end{equation}
 This last
inequality and (\ref{dim pro 1}) imply that $v''\not=v'.$

Since $\eta'_k \in H_{\eta'}$, the monotonicity of $T$ again gives
$$
\langle\xi_1(\eta'_k)-\xi_1(\eta'),v'_k-v'\rangle\ge 0, \qquad
\forall k\ge k_0;
$$
dividing by $\|\xi_1(\eta'_k)-\xi_1(\eta')\|$ and passing to the
limit, condition \emph{c.} in Lemma \ref{lemma} guarantees
\begin{equation}\label{dim pro 2bis}
\langle\xi_1(\eta)-\xi_1(\eta'),v''-v'\rangle\ge 0.
\end{equation}
Now summing the inequalities in (\ref{dim pro 2}) and in (\ref{dim
pro 2bis}), we obtain an inequality in contradiction with (\ref{dim
pro 1}). This concludes the proof.
 \QED

\section{Local boundedness of maximal H-monotone operators.}

It is well known that a maximal monotone operator
$T:\R^n\rightrightarrows\R^n$ is locally bounded. The proof relies essentially
on the fact that, given any ball, there exist $n+1$ points whose
convex hull contains the ball. In the case of operators
$T:\H^n\rightrightarrows V_1\sim\R^{2n}$ the situation is much more
involved. This section is essentially devoted to the proof Theorem
\ref{teo bounded}. We first show that a maximal H-monotone operator
defined on all $\Hn$ is locally bounded on every vertical segment
(see Proposition \ref{boundedness_on_vertical_segments}). We consider that this step
is really the bulk of the paper. Secondly, we show that  $T$
inherits the local boundedness on every horizontal segment from the
local boundedness of the vertical ones following an idea from \cite{BaRi2003}.

\begin{proposition}\label{boundedness_on_vertical_segments}
Let $T: \H^n \rightrightarrows V_1$ be a maximal $H$-monotone map
such that $\mathrm{dom}(T)=\H^n.$ Then the restriction of $T$ to any
vertical line is locally bounded, i.e. for every set of the type
$L:=\{\eta=(x,y,t)\in\H^n:\ t\in I\},$ with $x$ and $y$ fixed and
$I\subseteq \R,$ $I$ compact interval,  there exists $K=K(I)$ such
that $\diam_{\R^{2n}} \{T(L)\}\le K$.
\end{proposition}
\noindent{\bf Proof:} The proof is by contradiction. Assume that
there exists one of these vertical segments on which $T$ is not
bounded. Without loss of generality we can assume that the segment
is cointained in the $t$ axis. Moreover, we can assume that there exists a
sequence of points on the $t$ axis of the form $\eta_k=(0,0,h_k)$
such that $h_k \to 0$ and
\begin{equation} \label{lim}
 \lim_{k \to \infty} \diam_{\R^{2n}} \{ T(\eta_k) \} = \infty.
 \end{equation}
To obtain a contradiction we use a measure--theoretical argument as
follows.  Consider the sets:
$$ A_k = \{ \eta\in B_{\Hn} (e,5): \ \exists \  u \in T(\eta)\  \text{ such that} \ \|u\| > \frac{k}{2^{n}} \ \}.  $$
We will  construct measurable subsets $S_k\subset A_k$ with the property that there exists a constant $c>0$ such that for any $k$ we have
\begin{equation} \label{measure}
 \mathcal L^{2n+1} (S_k) >c,
 \end{equation}
where $ \mathcal L^{2n+1}$ is the Lebesgue measure in $\Hn$.

\medskip

Assuming the existence of $S_k$ let us show how to get the desired contradiction. We consider the sets
$$
U_k= \bigcup_{m\geq k} S_m.$$ Then $U_k$ is measurable and it is
decreasing: $U_{k+1} \subset U_k$ and $ \mathcal L^{2n+1} (U_k) >c$.
Let $S:= \bigcap_k U_k$, then $S$ is measurable and $ \mathcal
L^{2n+1} (S) \geq c$. Let $\eta\in S$, then $\eta$ lies in
infinitely many sets $S_k$. In particular there exists a sequence
$h(k)$ of indexes $h(k) \to \infty$ such that $\eta \in S_{h(k)}$
for each $h(k)$. This implies that there exists $u_{h(k)} \in
T(\eta)$ such that $ \|u_{h(k)}\| \geq h(k)$. On the other hand
$T(\eta)$ is compact by Proposition \ref{remark 1}, which is a
contradiction. \psn

In the following we will construct the sets $S_k\subset A_k$ (first
step) and will show the existence of a constant $c>0$ (independent
on $k$) for which \eqref{measure} holds for any $k$ (second step).

\medskip
\noindent \textbf{First step.}  The construction of $S_k$ uses the
measurable selection theorem (see e.g. \cite{Wag1977}). Let us
observe first that by \eqref{lim} it follows that $A_k \neq
\emptyset$ for any $k$. Moreover, for $k\geq 1$ there exists $h(k)
\in \mathbb N$ such that
$$ \diam_{\R^{2n}} \{ T(\eta_{h(k)}) \}  \geq 10 k.$$
To ease the notation we can assume that $h(k)=k$. We obtain a
sequence $\{u_k\}_k$ with $u_k \in T(\eta_k)$ such that
\begin{equation} \label{ten-n}
\|u_k\| \geq 10 k.
\end{equation}
Let us consider the unit vector in $V_1$
\begin{equation} \label{omega n}
\omega_k = \frac{u_k}{\|u_k\|},
\end{equation}
and the horizontal segment
$$
L_k= \{ \nu_k(t):=\eta_k\circ\exp(t\omega_k)\in\H^{n}:\ t \in [1,2] \}.
$$
We claim that $L_k \subset A_{10k}$. Indeed, let $\nu_k(t) \in L_k$
and $v_k \in T(\nu_k(t))$. By the $H$-monotonicity of $T$ we have
$$
\langle \xi_1(\nu_k(t))-\xi_1(\eta_k),v_k-u_k\rangle\ge 0
$$
and hence, by (\ref{ten-n}) and (\ref{omega n})
$$
\langle  v_k,t \omega_k\rangle \geq \langle u_k, t \omega_k\rangle \geq 10 k.
$$
Since $\omega_k$ is a unit vector by the Cauchy-Schwarz inequality
we obtain
\begin{equation} \label{v n}
\| v_k\| \geq \langle  v_k,\omega_k\rangle   \geq 10 k.
\end{equation}

The idea of the proof is to enlarge the segment $L_k$ by glueing
$2n$-dimensional sectors in the horizontal plane of each of its
points. We will prove that by this construction we obtain an
enlarged $(2n+1)$-dimensional set which is still a subset of $A_k$ and
whose Lebesgue measure is bounded below by a uniform constant.

\noindent Let us consider $I\subset \R^{2n-1}$ given by $$
I=[0,\pi]\times [0,\pi]\times \dots \times [0,\pi]\times [0,2\pi)
$$
and the spherical coordinates $\omega:I\to S^{2n-1}$ given by
$\omega(\Phi)=(\omega^1(\Phi),\ldots,\omega^{2n}(\Phi)),$
 for $\Phi=(\phi^1,\ldots,\phi^{2n-1}),$
\begin{equation}\label{def w}
\left\{\begin{array}{l}
\omega^1(\Phi)=\cos \phi^1\\
\omega^2(\Phi)=\sin \phi^1\cos\phi^2\\
\dots  \dots\\
\omega^{2n-1}(\Phi)=\sin \phi^1\sin \phi^2\dots \sin \phi^{2n-2}\cos\phi^{2n-1}\\
\omega^{2n}(\Phi)=\sin \phi^1\sin \phi^2\dots \sin \phi^{2n-2}\sin\phi^{2n-1}\\
\end{array}\right.
\end{equation}

To carry out the proposed construction let us select for each
$t\in [1,2]$ a vector $\tilde v_k(t) \in T(\nu_k(t))$. Here we
apply the measurable selection theorem (see \cite{Wag1977}) to
obtain, for every $k$, a measurable map $t\to \tilde v_k(t)$. In
the following consideration we will fix the index $k$. However
we note that, by (\ref{v n}),
 \begin{equation}\label{tilde theta0}
\| \tilde v_k(t)\| \geq 10 k.
\end{equation}
For each $t\in [1,2]$ let us write
\begin{equation}\label{tilde_theta}
\tilde v_k(t)=\|\tilde v_k(t)\|
\omega(\tilde\Phi_k(t)),
\end{equation}
for a suitable
$\tilde\Phi_k(t)=(\tilde\phi^1_k(t),\tilde\phi^2_k(t),\dots,\tilde\phi^{2n-1}_k(t))\in
I.$ Since the mapping $t\to \tilde v_k(t)$ is measurable we
obtain that the function
$
t \to (\tilde\Phi_k(t))
$
is measurable as
well. Set $\underline{i}=(i_1,i_2,\dots,i_{2n-1}),$ where $i_j\in \{1,2,3,4\}$ if $j\in \{1,2,\dots, 2n-2\},$ and $i_{2n-1}\in \{1,2,\dots, 8\},$ and
denote by $I^{\underline{i}}$
the set
$$
I^{\underline{i}}=\left[(i_1-1)\frac{\pi}{4},i_1\frac{\pi}{4}\right)\times
\left[(i_2-1)\frac{\pi}{4},i_2\frac{\pi}{4}\right)\times\dots\times
\left[(i_{2n-1}-1)\frac{\pi}{4},i_{2n-1}\frac{\pi}{4}\right).
$$
Fix any $\underline{i}$ as above, and consider the set
$$
T^{\underline{i}}_k= \left\{ t\in [1,2]: \tilde\Phi_k(t)\in I^{\underline{i}}\right\}.
$$
Then $\{T^{\underline{i}}_k\},$ for $k$ fixed, are disjoint
measurable sets with the property that
$\bigcup_{\underline{i}}T^{\underline{i}}_k= [1,2]$ up to a set of
measure 0. This implies that there exists $\underline{i}_0(k)$ such
that
$$
\mathcal L^{1}\left(T^{\underline{i}_{0}(k)}_k\right)\geq \frac{1}{2\cdot 4^{2n-1}}.
$$
Let us consider the subset of $L_{k}$ defined by
 \begin{equation}\label{nu n}
L^{\underline{i}_{0}(k)}_k = \{\nu_k(t):=\eta_k\circ\exp(t\omega_{k})\in\H^n:\ t
\in T^{\underline{i}_{0}(k)}_k \},
 \end{equation}
and to each $\underline{i}$ we associate the sector
$$
S^{\underline{i}}= \{\rho \omega(\Phi): \
\rho\in  [1,2], \ \Phi=(\phi^1,\phi^2,\dots,\phi^{2n-1})\in
I^{\underline{i}}\}.
$$
These sectors are $2n$-dimensional and disjoint. We define the
desired set $S_k$ by
\begin{equation} \label{addsectors}
 S_{k}= \{ \nu:=\nu_k(t)\circ\exp(\rho \omega(\Phi))\in \H^n: \ t\in T^{\underline{i}_0(k)}_k, \ \rho\omega(\Phi) \in S^{\underline{i}_0(k)} \}.
 \end{equation}
 It is clear that for $k$ sufficiently large, by the
 construction, we have $S_{k} \subset B_{\Hn}(e,5)$. We claim first
 that $S_{k}\subseteq A_{k}$. To see this let $\nu =
 \nu_k(t)\circ\exp(\rho \omega(\Phi))$ be an arbitrary point in
 $S_{k},$ and let $v \in T(\nu)$.

 We intend to prove that
 $\|v\|\geq \frac{k}{2^{n}}$.

  This will be done, using the fact that $\nu\in H_{\nu_k(t)}$ and the monotonicity of $T$ by comparing $(\nu,v)$ to
 the point
 $(\nu_k(t), \tilde v_k(t)),$ i.e.
$$
\langle \xi_1(\nu)-\xi_1(\nu_k(t)),v-\tilde v_k(t)\rangle\ge 0
$$
which implies
\begin{equation}\label{sti}
\rho\langle \omega(\Phi),v-\tilde
v_k(t)\rangle\ge 0
\end{equation}
 Let us note first that, if $\Phi=(\phi^1, \phi^2,\dots,
 \phi^{2n-1})$ and $\Psi=(\psi^1, \psi^2,\dots, \psi^{2n-1})$
 belong to the same $(2n-1)$-cube $I^{\underline{i}}$, then
\begin{equation}\label{estimate}
\langle \omega(\Phi), \omega(\Psi)\rangle \ge 2^{-(2n-1)/2}.
\end{equation}
Indeed, from the expression of the left hand side of the
previous inequality and taking into account the equation for the
the spherical coordinates, we have
\begin{align*}
\sum_{i=1}^{2n} \omega^i(\Phi)\omega^i(\Psi)=
&\cos\phi^1\cos \psi^1+\\
&+\sin \phi^1\cos\phi^2\sin \psi^1\cos \psi^2+\\
& +\dots +\\
&+\sin\phi^1\sin\phi^2\dots\sin \phi^{2n-2}\cos\phi^{2n-1} \sin\psi^1\sin\psi^2\dots\sin \psi^{2n-2}\cos\psi^{2n-1}+\\
&+\sin\phi^1\sin\phi^2\dots\sin \phi^{2n-2}\sin\phi^{2n-1} \sin\psi^1\sin\psi^2\dots\sin \psi^{2n-2}\sin\psi^{2n-1};
\end{align*}
if we take the last two lines of the sum above we have
\begin{align*}
\omega^{2n-1}(\Phi)\omega^{2n-1}(\Psi)&+\omega^{2n}(\Phi)\omega^{2n}(\Psi)=\\
&=\sin\phi^1\sin\phi^2\dots\sin \phi^{2n-2}\sin\psi^1\sin\psi^2\dots\sin \psi^{2n-2}\cos(\phi^{2n-1}-\psi^{2n-1})\\
& \ge \frac{1}{\sqrt 2}\sin\phi^1\sin\phi^2\dots\sin \phi^{2n-2}\sin\psi^1\sin\psi^2\dots\sin \psi^{2n-2},
\end{align*}
noticing that to obtain the previous inequality we use the fact that
$\sin \psi^i$ and $\sin \phi^i$ are nonnegative. Iterating this
argument, we finally get \eqref{estimate}.

\noindent Hence, by (\ref{tilde theta0}), (\ref{tilde_theta}) and
(\ref{sti}), and recalling that by definition of the set $S_k$ in
(\ref{addsectors}) we have that $\tilde\Phi_{k}(t)$ and $\Phi$ lie
in the same $(2n-1)$-cube $I^{\underline{i}_0(k)},$
 \begin{eqnarray*}
  \|v\|& \geq &\langle v, \omega(\Phi)\rangle\\
  & \geq & \langle \tilde v_k(t), \omega(\Phi)\rangle\\
&=&  \|\tilde v_k(t)\|\langle \omega(\tilde\Phi_k(t)), \omega(\Phi)
\rangle\\
&\ge &  10k {2}^{-(2n-1)/2}>\frac{k}{2^{n}}.
  \end{eqnarray*}

\medskip
  \noindent \textbf{Second step.}
  Our second claim is, that there exists a constant $c>0$ with the property that
  $$ \mathcal L^{2n+1}(S_{k}) \geq c.$$
  To prove this fact let us consider, for every $k,$ the mapping
  $$ F_k=(F^1_k,\ldots,F^{2n+1}_k): [1,2]\times [1,2]\times I^{\underline{i}_0(k)} \to \H^{n}, $$
  given by
  $$ F_k(t,\rho, \Theta) = \nu_k(t)\circ\exp (\rho \omega (\Theta)),
  $$
where $\nu_k(t)$ is as in (\ref{nu n}) and $\Theta=(\theta^1,\dots,
\theta^{2n-1}).$ Let $\Phi_k\in I$ be such that
$\omega_k=\omega(\Phi_k).$

Our aim is to show that if $\Theta$ is suitably chosen with respect
to $\Phi_k,$ then $|\texttt{\rm det}(JF_k(t,\rho,
 \Theta))|$ is bounded from below by a positive constant, where
$JF$ is the Jacobian of the function $F_k$.

\noindent Since $\omega_k$ is fixed, we can assume, without loss of
generality and to simplify the computations, that
$\omega_k=(1,0,\dots, 0),$ i.e. $\Phi_k=(0,\ldots,0)$.

\noindent Recalling that $\eta_k=(0,0,h_k),$ we obtain the
formula
$$ F_k(t, \rho, \Theta)=(F^1_k,\ldots,F^{2n+1}_k)(t, \rho, \Theta)=
\left(
\begin{array}{c}
 t+\rho \cos \theta^1\\
 \rho \sin \theta^1\cos \theta^2\\
 \rho \sin \theta^1\sin \theta^2 \sin \theta^3\cos \theta^4\\
 \dots\\
 \rho \sin \theta^1\sin \theta^2\sin \theta^3\dots \sin \theta^{2n-1}\\
 \rho \sin \theta^1\sin \theta^2\sin \theta^3\dots \cos \theta^{2n-1}\\
 h_k-2t\rho \sin \theta^1\sin \theta^2\dots \sin \theta^{n}\cos\theta^{n+1}
\end{array}\right)
$$
 Let us consider the Jacobian $JF_k$ of the function $F_k.$ If
 $n=1,$ trivial computations show that $|\texttt{\rm
 det}(JF_k(t,\rho, \theta))|=2\rho^2 |\sin\theta|.$ In the general
 case, we note that the first three columns of $JF_k$ are
 $$
 \left(
 \begin{array}{c}
 1\\
 0\\
 0\\
 \dots\\
0\\
0\\
-2\rho \prod_{i=1}^{n}\sin \theta^i\cos \theta^{n+1}
\end{array}\right),\quad
\left(
 \begin{array}{c}
 \cos\theta^1\\
 \sin\theta^1\cos\theta^2\\
 \sin\theta^1\sin\theta^2\cos\theta^3\\
 \dots\\
\prod_{i=1}^{2n-2}\sin \theta^i\cos\theta^{2n-1}\\
\prod_{i=1}^{2n-2}\sin \theta^i\sin\theta^{2n-1}\\
-2t \prod_{i=1}^{n}\sin \theta^i\cos \theta^{n+1}
\end{array}\right),\quad
\left(
 \begin{array}{c}
-\rho \sin\theta^1\\
 \rho\cos\theta^1\cos\theta^2\\
 \rho\cos\theta^1\sin\theta^2\cos\theta^3\\
 \dots\\
\rho\cos\theta^1\prod_{i=2}^{2n-2}\sin \theta^i\cos\theta^{2n-1}\\
\rho\cos\theta^1\prod_{i=2}^{2n-2}\sin \theta^i\sin\theta^{2n-1}\\
-2t\rho\cos\theta^1 \prod_{i=2}^{n}\sin \theta^i\cos \theta^{n+1}
\end{array}\right)
$$
In particular, the second and the third ones can be written as
$$
\sin\theta^1\left(
 \begin{array}{c}
 \cos\theta^1/\sin\theta^1\\
\cos\theta^2\\
\sin\theta^2\cos\theta^3\\
 \dots\\
\prod_{i=2}^{2n-2}\sin \theta^i\cos\theta^{2n-1}\\
\prod_{i=2}^{2n-2}\sin \theta^i\sin\theta^{2n-1}\\
-2t \prod_{i=2}^{n}\sin \theta^i\cos \theta^{n+1}
\end{array}\right),\quad\quad
\frac{\rho}{\cos\theta^1 }\left(
 \begin{array}{c}
-\sin\theta^1/\cos\theta^1\\
 \cos\theta^2\\
\sin\theta^2\cos\theta^3\\
 \dots\\
\prod_{i=2}^{2n-2}\sin \theta^i\cos\theta^{2n-1}\\
\prod_{i=2}^{2n-2}\sin \theta^i\sin\theta^{2n-1}\\
-2t\prod_{i=2}^{n}\sin \theta^i\cos \theta^{n+1}
\end{array}\right),
$$
therefore, if we remove the first entry, we get two dependent
columns. This means that, when computing the determinant of $JF_k$
starting from the first column, we have actually only one term to
consider, namely
$$
\texttt{\rm det}(JF_k(t,\rho,
 \Theta))=
-2\rho \prod_{i=1}^{n}\sin \theta^i\cos \theta^{n+1}\cdot \texttt{\rm det}(Jw^\rho(\rho,\Theta)),
$$
where $w^\rho(\rho,\Theta)=\rho w(\Theta)$ denotes the
$2n$-dimensional spherical coordinates (see (\ref{def w})). By
known computations,
$$
\texttt{\rm det}(Jw^\rho(\rho,\Theta))= \rho^{2n-1}\sin^{2n-2}\theta^1\sin^{2n-3}\theta^2\cdots\sin \theta^{2n-2};
$$
thus $$ |\texttt{\rm det}(JF_k(t,\rho,\Theta))|= 2\rho^{2n}
|\prod_{i=1}^{n}\sin \theta^i\cos
\theta^{n+1}\sin^{2n-2}\theta^1\sin^{2n-3}\theta^2\cdots\sin
\theta^{2n-2}|.
$$
We note that $\texttt{\rm det}(JF_k(t,\rho,\Theta))\not=0$
for a.e. $\Theta\in I^{\underline{i}_0(k)}$.
Let us consider the set
$$
C_k=T^{\underline{i}_0(k)}_k\times [1,2]\times
I^{\underline{i}_0(k)}.
$$
Since $S_k=F([1,2]\times [1,2]\times I^{\underline{i}_0(k)}),$ we
have that $F(C_k)\subseteq S_k.$ By the change of variable
formula we have that
$$
\mathcal L^{2n+1}(S_{k}) \geq \int_{C_{k}} |\texttt{\rm det}(JF_k(t,
\rho, \Theta))| dt d\rho d\theta^1 d\theta^2 \cdots d\theta^{2n-1}
\geq c_{0}\mathcal L^{2n+1}(C_{k}) := c. $$
 It is an exercise to show that $c$ is a uniform constant which does not depend on $k$, finishing the proof.
\QED

We are now able to prove Theorem \ref{teo bounded}:

\noindent{\textbf{Proof of Theorem \ref{teo bounded}}: We show
that $T$ is bounded in a suitable neighbourhood of the origin.
Let us consider the $4n$ segments in $\H^n:$
\begin{eqnarray*}
&&I^+_j:=\{(e_j,0,s)\in\H^n:\ -1\le s\le 1\}, \quad I^-_{j}:=\{(-e_j,0,s)\in\H^n:\
-1\le s\le 1\},\\
&& I^+_{j+n}:=\{(0,e_j,s)\in\H^n:\ -1\le s\le 1\},\quad    I^-_{j+n}:=\{(0,-e_j,s)\in\H^n:\
-1\le s\le 1\},
\end{eqnarray*}
where $j=1,2,\dots, n.$ Here $e_j$ denotes the $n$-tuple with $1$ in the $j$ position, and $0$ otherwise.

From Proposition \ref{boundedness_on_vertical_segments}, there is
$K>0$ such that $T(I^+_j), T(I^-_j)\subseteq B_{\R^{2n}}(0,K),$ for
every $j=1,\dots, 2n.$ Let $r\in (0,1)$ small enough such that, for
every $\xi\in B_{\H^n}(0,r)$ and for every $j=1,\dots, 2n,$ we have
$H_\xi\cap I^+_j, H_\xi\cap I^-_j\neq \emptyset:$ we note that by
the continuity of the map $\xi\mapsto H_\xi$ such $r>0$ exists since
the claim holds for $\xi=0$. Now, for any $\xi=(x,y,t)\in
B_{\H^n}(0,r)$ we define $\xi^+,\xi^-,\ v_j^+$ and $v_j^-$ by
$$
\xi^+_j=\xi\circ \exp(v^+_j(\xi))=H_\xi\cap I^+_j, \quad \xi^-_j=\xi\circ \exp(v^-_j(\xi))=H_\xi\cap I^-_j, \quad j=1,2,\dots, 2n.
$$
Straightforward computations show that $v^\pm_j$ coincide with one of the  vectors from the following list
$$
(e_j-x,-y),\quad (-e_j-x, -y),\quad (-x, e_j-y), \quad (-x,-e_j-y),
$$
thus $\|v^\pm_j\|\le 2,$ for every $j,$ and for every $\xi\in B_{\H^n}(0,r).$

From the H-monotonicity of $T,$ we have that
\begin{equation}\label{inequalities}
\langle u,v^\pm_j(\xi)\rangle \le  \langle u_j,v^\pm_j(\xi)\rangle\le K
\|v^\pm_j(\xi)\|\le 2K,
\end{equation}
for every $u\in T(\xi),$ $u_j\in T(\xi_j),$ and for every $j=1,\dots, 2n.$ The
inequalities \eqref{inequalities} imply that $T(\xi)$ is contained
in the polyhedron $P(\xi)$ defined by:
$$
P(\xi):=\{u\in V_1:\;\langle u,v^+_j(\xi)\rangle\le 2K, \langle u,v^-_j(\xi)\rangle\le 2K \quad
j=1,\dots, 2n\}.
$$
Note that there is no $v\in \R^{2n}\setminus \{0\}$ such that the half-space $\{u\in \R^{2n}:\, \langle v,u\rangle \le
0\}$ contains all the
vectors $\{v^\pm_j\}_{j=1,\dots,2n};$ as a consequence, the set $P(\xi)$ turns out to be a
polytope, i.e. it is bounded. Indeed, on the contrary, if $v\in
\R^{2n}\setminus \{0\}$ is such that $tv\in P(\xi),$ for every $t\ge
0,$ then $\langle v,v^\pm_j(\xi)\rangle \le 0,$ i.e. the set
$\{v^\pm_j(\xi)\}_j$ belongs to the half-space $\{u:\, \langle
v,u\rangle\le 0\},$ a contradiction. The continuity of the maps
$\xi\mapsto v^\pm_j(\xi),$ for every $j,$ entails, in particular, that
the set-valued map $\xi\mapsto P(\xi)$ is upper semicontinuous; thus,
if $r$ is small enough, there exists $K'\ge 2K$ such that
$$
P(\xi)\subseteq B_{\R^{2n}}(0,K'),\qquad \forall \xi\in B_{\H^n}(0,r).
$$
This implies that $T(\xi)\subseteq B_{\R^{2n}}(0,K'),$ for all $\xi\in
B_{\H^n}(0,r),$ therefore $T$ is locally bounded at the origin.
 \QED

Clearly, Theorem \ref{teo usc} and Theorem \ref{teo bounded}
give

\begin{corollary}\label{coroll}
Let $T: \H^n \rightrightarrows V_1$ be a maximal $H$-monotone
map, such that $\mathrm{dom}(T)=\H^n.$ Then $T$ is locally
bounded and upper semicontinuous.
\end{corollary}

\section{On Minty's theorem.}
This section we apply our main result in Theorem \ref{teo bounded}
in order to prove a horizontal version of Minty's theorem.  In the
following, for a given operator $T:\H^{n}\rightrightarrows V_1$ and
$\lambda>0,$ we denote by $T_\lambda:\H^n\rightrightarrows V_1$ the
operator
$$
T_\lambda=\xi_1+\lambda T.
$$
 It is clear that if $T$ is H-monotone, then $T_\lambda$ is
strictly H--monotone. We recall that in \cite{CaPi2014} the authors
prove Theorem \ref{Minty 1}, a result of Minty type in the case
$n=1$. Now, our aim is to prove Theorem \ref{Minty Hn}. In comparison
to Theorem \ref{Minty 1}, we will remove the
H-cyclically monotone assumption in i., and we also show that the result
holds for $\H^n.$ We note that in the Example
 \ref{example1} we have a map that satisfies the assumption in
Theorem \ref{Minty Hn}, but not the assumption i. in Theorem
\ref{Minty 1}.

  In order to prove the following result, we will follow the idea in \cite{BaCaKr2013}
  by using
degree-theoretical arguments for set valued maps \cite{HuPa1995}; the
results needed in the proof are collected in the Appendix of
\cite{BaCaKr2013}.

\medskip
\noindent \textbf{Proof of Theorem \ref{Minty Hn}}:

\noindent Let us first prove that $i.$ implies $ii,$ which is the more difficult part. Let $T$ be
maximal H-monotone with dom$(T)=\Hn$. Let us fix $\eta\in \H^n$ and
$\lambda>0.$ We consider the linear projection map $\pi:H_\eta\to V_1=\R^{2n}$
defined by $\pi(x,y,t)=(x,y).$ Note that since we restricted the projection to a hyperplane we
have that $\pi$ is bijective and we denote by $\pi^{-1}:\R^{2n}\to H_{\eta}$ its inverse.
We introduce the following notations:
$\widetilde {T_\lambda}$ is the operator $\widetilde
{T_\lambda}=T_\lambda\circ\pi^{-1}:\R^{2n}\rightrightarrows V_1$ and
$\pi(\zeta)=\tilde \zeta,\ \forall \zeta\in H_\eta.$ We have to
prove that $\widetilde {T_\lambda}$ is surjective.

Let us fix $p_0\in V_1\cong \R^{2n}:$ we show that it is possible to
find $R_0>0$ large enough such that
\begin{equation}\label{dim 0}
p_0\in \widetilde {T_\lambda}\left(B_{\R^{2n}}(\widetilde
\eta,R_0)\right):
\end{equation}
 in particular, we show this for
\begin{equation}\label{def R0}
R_0> \|p_0\|+\|\xi_1(\eta)\|+\lambda\sup\{\|v_\eta\|:\ v_\eta\in T(\eta)\}.
\end{equation}

Note, that the fact that the expression on the right in the above inequality  is finite follows from local boundedness of $T$.

\noindent{\bf Step 1.} In order to prove (\ref{dim 0}), we show first that
\begin{equation}\label{dim 1}\texttt{\rm deg}_{SV}\left(\widetilde {T_\lambda}-p_0,B_{\R^{2n}}(\widetilde \eta,R_0),0\right)=1,
\end{equation}
where $\texttt{\rm deg}_{SV}$ denotes the degree function for
set-valued maps. We consider the parametric set-valued map $\mathcal
F:[0,1]\times\overline {B_{\R^{2n}}(\widetilde
\eta,R_0)}\rightrightarrows V_1$  defined by
$$\mathcal F(\alpha, \tilde \zeta)=\tilde \zeta-p_0+\lambda\alpha T(\pi^{-1}(\tilde \zeta)),
$$ for all $\alpha\in [0,1],\ \tilde \zeta\in
\overline{B_{\R^{2n}}(\widetilde \eta,R_0)}$.

\noindent First we note that, by Proposition \ref{remark 1},  the
map $\mathcal F$ is convex-valued and compact-valued, i.e. for every
fixed $(\alpha ,\tilde\zeta)\in [0,1]\times
B_{\R^{2n}}(\widetilde\eta,R_0)$, the set $\mathcal F(\alpha,
\tilde\zeta)$ is compact and convex in
  $\mathbb R^{2n}$.
Moreover, Corollary \ref{cor_usc_closed} and Corollary \ref{coroll}
imply that $$\overline{\{\cup \mathcal F(\alpha, \tilde\zeta):\
(\alpha,\tilde\zeta)\in [0,1]\times \overline{B_{\R^{2n}}(\widetilde
\eta,R_0)}\}}$$ is compact in $\mathbb R^{2n}$. Finally, Corollary
\ref{cor_usc_closed} implies that the map
$(\alpha,\tilde\zeta)\mapsto \mathcal F(\alpha, \tilde \zeta)$ is
usc from $[0,1]\times \overline{B_{\R^{2n}}(\widetilde \eta,R_0)}$
into $2^{\mathbb R^{2n}}\setminus
  \{\emptyset\}$.

\noindent Now we are in the position to apply the mentioned
degree-theoretical arguments for set valued maps. According to the
above discussion, it follows that our map $\mathcal F(\alpha,\cdot)$ is a homotopy of class (P)
(see \cite{HuPa1995} and also Appendix in \cite{BaCaKr2013}).  The argument is based on the
application of Theorem 6.2 in \cite{BaCaKr2013}. In order to apply this statement  we need to show that the
constant curve $\gamma:[0,1]\to \R^{2n},$ defined by
$\gamma(\alpha)=0,$ is such that
\begin{equation}\label{deg 1}
\gamma(\alpha)\not\in \mathcal F(\alpha, \partial
B_{\R^{2n}}(\widetilde \eta,R_0)),\qquad \forall\alpha \in[0,1].
\end{equation}
We show (\ref{deg 1}) through arguing by contradiction: suppose that for some $\alpha$ there exists $\tilde\zeta\in
\partial B_{\R^{2n}}(\widetilde \eta,R_0)$ such that $$0\in \mathcal
F(\alpha, \tilde\zeta)=\tilde\zeta-p_0+\lambda\alpha
T(\pi^{-1}(\tilde\zeta)),$$ i.e. $p_0=\xi_1(\zeta)+\lambda\alpha
w_\zeta$ for some $w_\zeta\in T(\zeta),\ \zeta\in H_\eta$ and
$\zeta\in
\partial B_{\H^n}(\eta,R_0).$ This implies that, for every $v_\eta\in T(\eta),$
we have
$$p_0-\xi_1(\eta)-\lambda\alpha v_\eta=\xi_1(\zeta)-\xi_1(\eta)+\lambda\alpha( w_\zeta-v_\eta).$$
Multiplying the previous vector equality by
$(\xi_1(\zeta)-\xi_1(\eta))$ we obtain
$$\langle \xi_1(\zeta)-\xi_1(\eta),p_0-\xi_1(\eta)-\lambda\alpha v_\eta\rangle=\|\xi_1(\zeta)-\xi_1(\eta)\|^2+\lambda\alpha\langle
\xi_1(\zeta)-\xi_1(\eta),w_\zeta-v_\eta\rangle.$$
The H--monotonicity of $T$ implies
$$\|\xi_1(\zeta)-\xi_1(\eta)\|^2\le\left\|
\langle \xi_1(\zeta)-\xi_1(\eta),p_0-\xi_1(\eta)-\lambda\alpha
v_\eta\rangle\right\|\le\left\|
\xi_1(\zeta)-\xi_1(\eta)\right\|\cdot\left\|p_0-\xi_1(\eta)-\lambda\alpha
v_\eta\right\|;$$ hence $$ R_0\le
\left\|p_0-\xi_1(\eta)-\lambda\alpha v_\eta\right\|$$
 This contradicts (\ref{def R0}) and hence (\ref{deg 1}) holds.
The homotopy invariance property for $\mathcal F$ (see
Theorem 6.2 in \cite{BaCaKr2013}) gives that
$$
\texttt{\rm deg}_{SV}\left( \mathcal F(\alpha,
\cdot),B_{\R^{2n}}(\widetilde \eta,R_0),\gamma(\alpha) \right)$$ does
not depend on $\alpha:$ hence,
\begin{eqnarray}
\texttt{\rm deg}_{SV}\left(
\widetilde{T_\lambda}-p_0,B_{\R^{2n}}(\widetilde \eta,R_0),0 \right)
&=&\texttt{\rm deg}_{SV}\left( \mathcal F(1,
\cdot),B_{\R^{2n}}(\widetilde \eta,R_0),0 \right)\nonumber\\
&=&\texttt{\rm deg}_{SV}\left( \mathcal F(0,
\cdot),B_{\R^{2n}}(\widetilde \eta,R_0),0 \right)\nonumber\\
&=&\texttt{\rm deg}_{SV}\left( I_{\R^{2n}}-p_0,B_{\R^{2n}}(\widetilde\eta,R_0),0 \right)\nonumber\\
&=&\texttt{\rm deg}_{B}\left( I_{\R^{2n}}-p_0,B_{\R^{2n}}(\widetilde\eta,R_0),0 \right)\nonumber\\
&=&\texttt{\rm deg}_{B}\left(
I_{\R^{2n}},B_{\R^{2n}}(\widetilde\eta,R_0),p_0 \right),\label{dim 3}
\end{eqnarray}
where $\texttt{\rm deg}_{B}$ denotes the degree function for single-valued maps. Note that (\ref{def R0}) implies that $p_0\in
B_{\R^{2n}}(\widetilde\eta,R_0):$ hence $ \texttt{\rm deg}_{B}\left(
I_{\R^{2n}},B_{\R^{2n}}(\widetilde\eta,R_0),p_0 \right)=1$ (see Theorem
6.1 in \cite{BaCaKr2013}) and hence, by (\ref{dim 3}), we have that
(\ref{dim 1}) is true.

\noindent{\bf Step 2.}  By Step 1 and the definition of $\texttt{\rm
deg}_{SV}$, for small $\varepsilon>0$, one has that
\begin{equation}\label{deg-utolso}
    \texttt{\rm deg}_B(f_\varepsilon-p_0,B_{\R^{2n}}(\widetilde\eta,R_0),0)=1,
\end{equation}
where $f_\varepsilon: \overline{B_{\R^{2n}}(\widetilde\eta,R_0)}\to
\mathbb R^{2n}$ is a continuous approximate selector of the upper
semicontinuous set-valued map $\widetilde{T_\lambda}$
 such that
\begin{equation}\label{eps-approx}
    f_\varepsilon(\tilde\zeta)\in \widetilde{T_\lambda}\left(B_{\mathbb
R^{2n}}(\tilde\zeta,\varepsilon)\cap
\overline{B_{\R^{2n}}(\widetilde\eta,R_0)}\right)+B_{\mathbb
R^{2n}}(0,\varepsilon),\ \ \forall \tilde\zeta\in
\overline{B_{\R^{2n}}(\widetilde\eta,R_0)},
\end{equation}
see Proposition 6.1 in \cite{BaCaKr2013}.  Let
$\varepsilon=\frac{1}{k}$ and let $\phi_k:=f_{1/k}$, $k\in \mathbb
N.$ First of all, from (\ref{deg-utolso}) and the properties of the
Brouwer degree function ${\rm deg}_B$ (see Theorem 6.1 in \cite{BaCaKr2013}),
we have that for every $k\in \mathbb N$ there exists
$\tilde\zeta_k\in B_{\R^{2n}}(\widetilde\eta,R_0)$ such that
$p_0=\phi_k(\tilde\zeta_k)$. Up to a subsequence, we may assume that
$\tilde\zeta_k\to \tilde\nu\in
 \overline{B_{\R^{2n}}(\widetilde\eta,R_0)}$. On the
other hand, by relation (\ref{eps-approx}), we have that
$$
p_0=\phi_k(\tilde\zeta_k)\in
\widetilde{T_\lambda}\left(B_{\mathbb R^{2n}}(\tilde\zeta_k,1/k)\cap
\overline{B_{\R^{2n}}(\widetilde\eta,R_0)}\right)+B_{\mathbb
R^{2n}}(0,1/k) ,
$$
i.e., there exists $\tilde\nu_k\in B_{\mathbb
R^{2n}}(\tilde\zeta_k,\frac{1}{k})\cap
\overline{B_{\R^{2n}}(\widetilde\eta,R_0)}$ and $p_k\in B_{\mathbb
R^{2n}}(0,\frac{1}{k})$ such that $p_0\in
\widetilde{T_\lambda}(\tilde \nu_k)+p_k.$ Clearly, $\tilde\nu_k\to
\tilde\nu$ and $p_k\to 0$ as $k\to \infty$.
Let us consider the sequence $\{(p_0-p_k, \tilde\nu_k)\}\in \mathrm{graph}(\tilde{T_\lambda});$ by the usc of
$\tilde{T_\lambda}$ we get that $p_0\in \tilde{T_\lambda}(\tilde\nu).$

Finally, we claim that  $\tilde\nu\in B_{\R^{2n}}(\widetilde\eta,R_0).$ To see
this, let us assume, by contradiction, that  $\tilde\nu\in
\partial B_{\R^{2n}}(\widetilde\eta,R_0).$ Then, $p_0\in
\widetilde{T_\lambda}(\tilde\nu)$ is equivalent to $0\in
\widetilde{T_\lambda}(\tilde\nu)-p_0= \mathcal F(1, \tilde\nu)$,
which contradicts relation (\ref{deg 1}). Consequently,
$\tilde\nu\in B_{\R^{2n}}(\widetilde\eta,R_0)$; therefore we
obtain (\ref{dim 0}), which concludes the proof of the first
implication  $i.\ \Rightarrow\ ii$ of Theorem \ref{Minty Hn}.

 Let us prove that $ii.$ implies $i.$ The proof is essentially
in \cite{CaPi2014}, where the case $n=1$ is considered; however, for
the sake of completeness, we include it. Let $T:\Hn\rightrightarrows
V_1$ be a set-valued H-monotone map, with domain $\Hn,$ such that,
for every $\eta_0\in \Hn,$
$$
\mathrm{rge}(T_\lambda)|_{H_{\eta_0}}=V_1.
$$
We argue by contradiction and suppose that $T$ is not maximal
H-monotone. Then there exist $\eta_0\in \Hn,$ and $w\notin
T(\eta_0)$ such that, for every $\eta\in H_{\eta_0},$ and $v\in
T(\eta),$
\begin{equation}\label{inequality}
\langle w -v,\xi_1(\eta_0)-\xi_1(\eta)\rangle \ge 0.
\end{equation}
Without loss of generality, we assume that $\eta_0=e$: in fact, via
a left translation, the map $\eta\mapsto T(\eta_0\circ\eta)$ has the
same properties of $T.$ From the assumptions, for $\lambda=1$ we
have
 $\mathrm{rge}(T+\xi_1)|_{H_e}=V_1; $ therefore
\begin{equation}\label{equality}
w=\tilde v+\xi_1(\tilde \eta),
\end{equation}
for some $\tilde \eta\in H_e$ and $\tilde v\in T(\tilde \eta).$
 From \eqref{equality}, choosing $\eta=\tilde \eta$ in
 \eqref{inequality}, we
obtain
$$
-\langle \xi_1(\tilde \eta),\xi_1(\tilde \eta)\rangle \ge 0,
$$
i.e., $\xi_1(\tilde \eta)=0.$ Since $\tilde \eta\in H_e,$ we
deduce that $\tilde \eta=0,$ and $w=\tilde v\in T(e),$
contradicting our assumption on $w.$ This concludes the proof of
Theorem \ref{Minty Hn}.
 \QED

 \subsection{Lipschitz continuity of the resolvent operator in the Hausdorff metric. }\label{subsection resolvent}

In this subsection we are interested in studying the regularity of
the resolvent $Q_\lambda$ of a maximal H--monotone operator $T$
 defined by $$
Q_\lambda=(\xi_1+\lambda T)^{-1}:V_1\rightrightarrows\Hn.$$ First,
we have to recall that if $T$ is maximal H-monotone and
$\eta\in\Hn,$ then the map $T_\lambda|_{H_\eta}$ is not
injective, in general, and hence $\left(T_\lambda
|_{H_\eta}\right)^{-1}:V_1\rightrightarrows H_\eta$ is not
single-valued (see Example \ref{example2} below). Using the strictly
H--monotonicity of the operator $T_\lambda,$ the only information we have is that, for every $\eta'\in H_\eta,$ 
$$
T_\lambda(\eta)\cap T_\lambda(\eta')=\emptyset.
$$
We note that, for every fixed $v$, $Q_\lambda(v)$ is a closed subset of $\Hn,$
since it is the inverse image via the usc map $T_\lambda$ of a point.
Moreover, Theorem \ref{Minty Hn} implies that for every fixed $v\in
V_1$ and $\eta\in \Hn$, there exists at least one point $\eta'\in
H_\eta$ such that $v\in T_\lambda(\eta'),$ i.e. $\eta'\in
Q_\lambda(v)$. Therefore $H_{(0,0,h)}\cap
Q_\lambda(v)\not=\emptyset,$ for every $h\in\R.$ Hence
$Q_\lambda(v)$ is unbounded for every fixed $v\in V_1$. We summarize this discussion in the following:

\begin{remark}
 Let $T:\Hn\rightrightarrows V_1$ be a maximal H--monotone map
with $\dom(T)=\Hn.$ Then, for every $\lambda >0,$ the resolvent
$Q_\lambda:V_1\rightrightarrows \Hn$ is closed--valued, and $Q_\lambda(v)$ is unbounded for every $v\in V_1.$
\end{remark}

As we mentioned in the introduction, if we consider the resolvent in
our context, we are very far from the Euclidean situation where the
resolvent map $(I+\lambda T)^{-1}$ of a maximal monotone set-valued
map $T:\R^n\rightrightarrows \R^n$ is single-valued on $\R^n$ and
1-Lipschitz continuous. However, in this line of investigation, it
is useful to think about the
 notion of multivalued Lipschitz  map.

Let $Q:V_1\rightrightarrows\Hn$ be a closed--valued multivalued map.
We recall (see Definition 9.26 in \cite{RoWe2004}) that $Q$ is
\emph{Lipschitz continuous in the Hausdorff metric}, if
$\dom(Q)=V_1$ and there exists a positive $k$ such that
$$
Q(v')\subseteq Q(v)+B_{\Hn}(0,k\|v'-v\|),\quad \forall v,v'\in V_1.
$$
We have the following regularity result for our resolvent:

\begin{proposition}\label{Lip Q}
Let $T:\Hn\rightrightarrows V_1$ be a maximal H--monotone map with
$\dom(T)=\Hn.$ Then, for every $\lambda>0,$ the resolvent
$Q_\lambda$ is 1-Lipschitz continuous in the Hausdorff metric.
\end{proposition}
\textbf{Proof:} Let us consider $v$ and $v'$ in $V_1,$ with
$v\not=v'.$ For every $\eta\in Q_\lambda(v),$ i.e.
\begin{equation}\label{v}
v\in \xi_1(\eta)+\lambda T(\eta),
\end{equation}
Theorem \ref{Minty Hn} guarantees that there exists $\eta'\in
H_\eta$ such that $\eta'\in Q_\lambda(v'),$ i.e.
\begin{equation}\label{v'}
v'\in \xi_1(\eta')+\lambda T(\eta').
\end{equation}
Relations (\ref{v}) and (\ref{v'}) give that $ v- \xi_1(\eta)\in
\lambda T(\eta)$ and $ v'- \xi_1(\eta')\in \lambda T(\eta')$. Since
$\lambda T:\H\rightrightarrows V_1$ is H--monotone, we have
$$
\langle v- \xi_1(\eta)-v'+ \xi_1(\eta'),\xi_1(\eta)-
\xi_1(\eta')\rangle\ge 0
$$
and hence
 \begin{eqnarray*}
\|v-v'\|&\ge& \langle v- v', \frac{\xi_1(\eta)- \xi_1(\eta')}{\|
\xi_1(\eta)- \xi_1(\eta')\|}\rangle\\
&\ge& \| \xi_1(\eta)- \xi_1(\eta')\|.
\end{eqnarray*}
The previous inequality implies that for every $\eta\in
Q_\lambda(v)$ there exists $\eta'\in Q_\lambda(v')$ such that
$d_H(\eta,\eta')\le \|v-v'\|.$ This implies 1-Lipschitz continuity
of $Q_\lambda$  in the Hausdorff metric. The claim is proved.
 \QED

Let us conclude with the following example presented in
\cite{CaPi2014}:

\begin{example}\label{example2}
 Let us consider the gauge function $N:\H\to \R$ defined as
$$
N(x,y,t): =((x^2+y^2)^2+t^2)^{1/4}.
$$
It is known that this function is H-convex (see
\cite{DaGaNh2003}). The associated horizontal subgradient map
$\partial_HN$ is given by
$$
\partial_HN(x,y,t)=
     \begin{cases}
    \overline{B_{\R^2}(0,1)} & (x,y,t)=(0,0,0) \\
    \frac{1}{N^3(x,y,t)}\left(
    x(x^2+y^2)+yt,y(x^2+y^2)-xt\right) & (x,y,t)\neq (0,0,0).
    \end{cases}
    $$
 For every fixed $\lambda> 0,$ let the map
$T_\lambda:=\xi_1+\lambda\partial_H N:\H\rightrightarrows V_1$ that
is maximal strictly H--monotone. First, it is possible to prove that
there exist $\eta''\in \H,$ and $\eta, \eta'\in H_{\eta''},$
$\eta\neq \eta',$ such that
$$
T_\lambda (\eta)\cap T_\lambda(\eta')\neq \emptyset.
$$
Secondly, it is clear that $T_\lambda$ is not a Lipschitz
continuous map, i.e. it does not exist a positive $k$ such that
$$
T_\lambda(\eta')\subseteq
T_\lambda(\eta)+B_{\R^2}\left(0,k\,d_H(\eta,\eta')\right),\quad
\forall \eta,\eta'\in \H.
$$
In fact, for $\eta'=(0,0,0)$ and $\eta=(x,y,0)$ the previous
inclusion is false.

If we are interested in $Q_\lambda=(\xi_1+\lambda
 \partial_HN)^{-1}:V_1\rightrightarrows\H,$ an easy calculation gives $$Q_\lambda(0,0)=\{(0,0,t)\in\H;\
t\in\R\}.$$ Now, let us consider  $(x,y,t)\not=(0,0,0)$ and
$v\not=(0,0)$ with $(x,y,t)\in Q_\lambda(v),$ i.e.
$v=T_\lambda(x,y,t);$ straightforward computations lead to the
following
$$\epsilon^2\ge\|v\|^2=
\|T_\lambda(x,y,t)\|^2=
(x^2+y^2)\left(1+\frac{\lambda^2}{N^2(x,y,t)}+\frac{2\lambda}{N^3(x,y,t)}(x^2+y^2)\right);
$$
hence  $(x,y,t)\in Q_\lambda(v)$ implies $\|(x,y)\|\le \epsilon.$
Moreover, since $T_\lambda$ is surjective on every horizontal plane
$H_\eta$  and in particular on $H_{(0,0,t)},$  we obtain that for
every $t\not=0$ there exists $(x,y)$ such that $v\in
T_\lambda(x,y,t).$ These prove that $0\not=\|v\|\le \epsilon$ gives
$$Q_\lambda(v)\ \texttt{\rm is unbounded},\qquad
Q_\lambda\subset\{(x,y,t)\in\H:\ \|(x,y)\|\le \epsilon\}.
$$

\end{example}

\end{document}